\documentclass[12pt]{amsart}
\usepackage[margin=1.2in]{geometry}
\usepackage{graphicx,latexsym}
\usepackage{comment}
\usepackage{mathrsfs}
\usepackage{amsfonts, amssymb, amsmath, amsthm, amsxtra, amscd, hyperref}
\usepackage{cite}
\usepackage{mathtools}
\usepackage{todonotes}
\usepackage[foot]{amsaddr}
\usepackage{comment}
\usepackage{orcidlink}

\usepackage{tikz}
\makeatletter
\DeclareRobustCommand\widecheck[1]{{\mathpalette\@widecheck{#1}}}
\def\@widecheck#1#2{%
	\setbox\z@\hbox{\m@th$#1#2$}%
	\setbox\tw@\hbox{\m@th$#1%
		\widehat{%
			\vrule\@width\z@\@height\ht\z@
			\vrule\@height\z@\@width\wd\z@}$}%
	\dp\tw@-\ht\z@
	\@tempdima\ht\z@ \advance\@tempdima2\ht\tw@ \divide\@tempdima\thr@@
	\setbox\tw@\hbox{%
		\raise\@tempdima\hbox{\scalebox{1}[-1]{\lower\@tempdima\box
				\tw@}}}%
	{\ooalign{\box\tw@ \cr \box\z@}}}
\makeatother

\usetikzlibrary{arrows,chains,matrix,positioning,scopes}
\makeatletter
\tikzset{join/.code=\tikzset{after node path={%
			\ifx\tikzchainprevious\pgfutil@empty\else(\tikzchainprevious)%
			edge[every join]#1(\tikzchaincurrent)\fi}}}
\makeatother
\tikzset{>=stealth',every on chain/.append style={join},
	every join/.style={->}}
\tikzstyle{labeled}=[execute at begin node=$\scriptstyle,
execute at end node=$]
\allowdisplaybreaks

\newcommand{\N}{\mathbb{N}}

\newcommand{\R}{\mathbb{R}}
\newcommand{\C}{\mathbb{C}}

\newcommand{\dx}{{\rm d}x }

\newcommand{\supp}{\operatorname{supp}}

\newtheorem{theorem}{Theorem}[section]

\newtheorem{lemma}[theorem]{Lemma}
\newtheorem{corollary}[theorem]{Corollary}

\theoremstyle{definition}

\theoremstyle{remark}

\numberwithin{equation}{section}

\begin{document}
	
	\title[Linear topological invariants by shifted fundamental solutions]{Linear topological invariants for kernels of differential operators by shifted fundamental solutions}
	
	\author[A.~Debrouwere]{Andreas Debrouwere\orcidlink{0000-0003-4416-5758}}
	\address{%
		Department of Mathematics and Data Science\\
		Vrije Universiteit Brussel\\
		Pleinlaan 2\\
		1050 Brussels\\
		Belgium}
	
	\email{andreas.debrouwere@vub.be}

	\author[T.~Kalmes]{Thomas Kalmes\orcidlink{0000-0001-7542-1334}}
	\address{%
		Faculty of Mathematics\\
		Chemnitz University of Technology: Technische Universit\"at Chemnitz\\
		09107 Chemnitz\\
		Germany}
	
	\email{thomas.kalmes@math.tu-chemnitz.de}

	\begin{abstract} 
		
		We characterize the condition $(\Omega)$  for smooth kernels of partial differential operators in terms of the existence of shifted fundamental solutions satisfying certain properties. The conditions $(P\Omega)$ and $(P\overline{\overline{\Omega}})$  for distributional kernels are characterized in a similar way. By lifting theorems for Fr\'echet spaces and (PLS)-spaces, this provides characterizations of the problem of parameter dependence for smooth and distributional solutions of differential equations by shifted fundamental solutions.
		As an application, we give a new proof of the fact that  the space $\{ f \in \mathscr{E}(X) \, | \, P(D)f = 0\}$ satisfies $(\Omega)$ for any differential operator $P(D)$ and any open convex set  $X \subseteq \R^d$.\\

		\\ 
		
		
	\end{abstract}

	\maketitle
	
	\section{Introduction}
	
	In their seminal work \cite{MTV90}  Meise, Taylor and Vogt characterized the constant coefficient linear partial differential operators $P(D) = P(-i\frac{\partial}{\partial x_1},\ldots,-i\frac{\partial}{\partial x_d})$ that have a continuous linear right inverse on $\mathscr{E}(X)$ and/or $\mathscr{D}'(X)$ ($X \subseteq \R^d$ open) in terms of the existence of  certain shifted fundamental solutions of $P(D)$.  
	Later on, Frerick and Wengenroth \cite{FW04b, Wengenroth11}  gave a similar characterization of the  surjectivity of $P(D)$ on $\mathscr{E}(X)$, $\mathscr{D}'(X)$,  and $\mathscr{D}'(X)/\mathscr{E}(X)$ as well as  of the existence of right inverses of $P(D)$ on the latter space. Roughly speaking, these results assert that $P(D)$ satisfies some condition (e.g.\ being surjective on $\mathscr{E}(X)$) if and only if for each compact subset $K$ of $X$ and $\xi\in X$ far enough away from  $K$ there is a shifted fundamental solution $E$  for $\delta_\xi$ such that $E$ satisfies a certain property on $K$. Of course, this property depends on the condition one wants to characterize. Results of the same type have also been shown for spaces of non-quasianalytic ultradifferentiable functions and ultradistributions \cite{Langenbruch1996, Langenbruch1998} and for spaces of real analytic functions \cite{Langenbruch2004}. The aim of this paper is to complement the above results by characterizing several linear topological invariants for smooth  and distributional kernels of $P(D)$ by means of shifted fundamental solutions.
	
	The study of linear topological invariants for kernels of $P(D)$ goes back to the work of Petzsche \cite{Petzsche1980} and Vogt \cite{Vogt1983-2} and was reinitiated by Bonet and  Doma\'nski \cite{B-D2006,B-D2008,Domanski}. It is motivated by the question of surjectivity of $P(D)$ on  vector-valued function and distribution spaces, as we now proceed to explain. We assume that the reader is familiar with the condition $(\Omega)$ for Fr\'echet spaces \cite{M-V} and the conditions $(P\Omega)$ and $(P\overline{\overline{\Omega}})$ for (PLS)-spaces \cite{B-D2006,Domanski} (see also the preliminary  Section  \ref{sec: preliminaries}). Set $\mathscr{E}_P(X) = \{ f \in \mathscr{E}(X) \, | \, P(D)f = 0\}$ and  $\mathscr{D}'_P(X) = \{ f \in \mathscr{D}'(X) \, | \, P(D)f = 0\}$. Suppose that $P(D)$ is surjective on $\mathscr{E}(X)$, respectively, $\mathscr{D}'(X)$. Given a locally convex space $E$, it is natural to ask whether $P(D):\mathscr{E}(X;E)\rightarrow\mathscr{E}(X;E)$,
	respectively, $
	P(D):\mathscr{D}'(X;E)\rightarrow\mathscr{D}'(X;E)
	$
	is still surjective. If $E$ is a space of functions or distributions, this question is a reformulation of the well-studied problem of parameter dependence for solutions of partial differential equations; see \cite{B-D2006,B-D2008,Domanski} and the references therein.
	The splitting theory for Fr\'echet spaces \cite{Vogt1987} implies that the mapping $P(D):\mathscr{E}(X;E)\rightarrow\mathscr{E}(X;E)$ for $E =\mathscr{D}'(Y)$ ($Y \subseteq \R^n$ open) or $\mathscr{S}'(\R^n)$  is surjective if and only if $\mathscr{E}_P(X)$ satisfies $(\Omega)$. Similarly, as an application of their lifting results for (PLS)-spaces, Bonet and Doma\'nski showed that the mapping  $
	P(D):\mathscr{D}'(X;E)\rightarrow\mathscr{D}'(X;E)
	$ for $E =\mathscr{D}'(Y)$ or $\mathscr{S}'(\R^n)$  is surjective if and only if $\mathscr{D}_P(X)$ satisfies $(P\Omega)$ \cite{B-D2006}, while it is surjective for  $E =\mathscr{A}(Y)$  if and only if $\mathscr{D}'_P(X)$ satisfies $(P\overline{\overline{\Omega}})$ \cite{Domanski}.
	
	Petzsche \cite{Petzsche1980} showed that $\mathscr{E}_P(X)$ satisfies $(\Omega)$ for  any convex open set $X$\footnote{Petzsche actually showed this result under the additional hypothesis that $P(D)$ is hypoelliptic. However, as observed in \cite{DeKa22}, a careful inspection of his proof shows that this hypothesis can be omitted.}, while Vogt proved that this is the case for an arbitrary open set $X$ if $P(D)$ is elliptic \cite{Vogt1983-2}. Similarly,  $\mathscr{D}'_P(X)$ satisfies $(P\Omega)$ for any convex open set $X$  \cite{B-D2006} and for an arbitrary open set $X$ if $P(D)$ is elliptic \cite{B-D2006,FrKa,Vogt1983-2}. On the negative side, the second author \cite{Kalmes12-2} constructed a differential operator $P(D)$ and an open set $X\subseteq\R^d$ such that $P(D)$ is surjective on $\mathscr{D}'(X)$ (and thus also on $\mathscr{E}(X)$)  but $\mathscr{E}_P(X)$ and  $\mathscr{D}'_P(X)$ do not satisfy $(\Omega)$, respectively, $(P\Omega)$. Furthermore, $\mathscr{D}'_P(X)$ does not satisfy $(P\overline{\overline{\Omega}})$ for any convex open set $X$ if $P(D)$ is hypoelliptic  and for an arbitrary open set $X$ if $P(D)$ is elliptic \cite{Domanski,Vogt2006}. We refer to \cite{DeKa22, Domanski} and the references therein for further results concerning  $(\Omega)$ for $\mathscr{E}_P(X)$ and $(P\Omega)$ and $(P\overline{\overline{\Omega}})$ for $\mathscr{D}_P(X)$.
	
	Apart from this classical application to the problem of surjectivity of $P(D)$ on spaces of vector-valued smooth functions and distributions, in our recent article \cite{DeKa22-2}, the linear topological invariant $(\Omega)$ for $\mathscr{E}_P(X)$ played an important role to establish quantitative approximation results of Runge type for several classes of partial differential operators. See  \cite{EnPe21,RuSa19,RuSa20}
	for other works on this topic.

	In the present note, we characterize the condition $(\Omega)$ for $\mathscr{E}_P(X)$ and the conditions $(P\Omega)$ and $(P\overline{\overline{\Omega}})$ for $\mathscr{D}'_P(X)$ in terms of the existence of certain shifted fundamental solutions for $P(D)$.  By the above mentioned results from \cite{B-D2006,Domanski}, the latter provides characterizations of the problem of distributional and real analytic parameter dependence for distributional solutions of the equation $P(D) f = g$ by shifted fundamental solutions. This answers a question of Doma\'nski \cite[Problem 7.5]{Domanski2011} for distributions.

	We now state our main result. Set $\N=\{0,1,2,\ldots\}$. Let $Y  \subseteq \R^d$ be relatively compact and open. For $N\in \N$ we define
	$$
	\|f\|_{\overline{Y},N}= \max_{x\in \overline{Y},|\alpha|\leq N}|f^{(\alpha)}(x)|, \qquad f \in C^N(\overline{Y}),
	$$
	and
	$$
	\|f\|^*_{\overline{Y},N} =  \sup\{|\langle f,\varphi\rangle| \, | \,\varphi\in\mathscr{D}_{\overline{Y}}, \|\varphi\|_{\overline{Y},N}\leq 1\}, \qquad  f \in (\mathscr{D}_{\overline{Y}})',
	$$
	where $\mathscr{D}_{\overline{Y}}$ denotes the Fr\'echet space of smooth functions with support in  $\overline{Y}$. 
	
	\begin{theorem}\label{thm: main2}
		Let $P\in\C[\xi_1,\ldots,\xi_d]$, let $X\subseteq\R^d$ be open, and let $(X_n)_{n\in\N}$ be  an exhaustion by relatively compact open subsets of $X$. 
		\begin{itemize}
			\item[(a)]
			$P(D): \mathscr{E}(X) \to \mathscr{E}(X)$ is surjective and  $\mathscr{E}_P(X)$ satisfies $(\Omega)$ if and only if
			\begin{eqnarray}\label{eq: Omega}
				\begin{gathered}
					\forall\,n\in\N\,\exists\,m\geq n, N\in\N\,\forall\,k\geq m, \xi\in X\backslash \overline{X_m}\,
					\exists\,K\in\N, s,C>0\\
					\forall \,\varepsilon \in (0,1) \, \exists\,E_{\xi,\varepsilon}\in\mathscr{D}'(\R^d) \mbox{ with }  P(D)E_{\xi,\varepsilon}=\delta_\xi\mbox{ in }X_k  \\ 
					\mbox{ such that }  \|E_{\xi,\varepsilon}\|_{\overline{X_n},N}^*\leq \varepsilon \quad \mbox{and} \quad \|E_{\xi,\varepsilon}\|_{\overline{X_k},K}^*\leq \frac{C}{\varepsilon^s}.
				\end{gathered}
			\end{eqnarray}
			\item[(b)]	 $P(D): \mathscr{D}'(X) \to \mathscr{D}'(X)$ is surjective and  $\mathscr{D}'_P(X)$ satisfies $(P\Omega)$ if and only if
			\begin{eqnarray}\label{eq: POmega}
				\begin{gathered}
					\forall\,n\in\N\,\exists\,m\geq n\,\forall\,k\geq m, N\in\N, \xi\in X\backslash \overline{X_m} \, \exists\,K\in\N, s,C>0\\
					\forall\,\varepsilon \in (0,1) \, \exists\,E_{\xi,\varepsilon}\in\mathscr{D}'(\R^d)\cap C^N(\overline{X_n}) \mbox{ with }  P(D)E_{\xi,\varepsilon}=\delta_\xi\mbox{ in }X_k   \\
					\mbox{ such that }\|E_{\xi,\varepsilon}\|_{\overline{X_n},N}\leq \varepsilon \quad \mbox{and} \quad \|E_{\xi,\varepsilon}\|_{\overline{X_k},K}^*\leq \frac{C}{\varepsilon^s}.
				\end{gathered}
			\end{eqnarray}
			\item[(c)]	 $P(D): \mathscr{D}'(X) \to \mathscr{D}'(X)$ is surjective and  $\mathscr{D}'_P(X)$ satisfies $(P\overline{\overline{\Omega}})$  if and only if \eqref{eq: POmega} with $``\exists\,K\in\N, s,C>0''$ replaced by  $``\forall s > 0 \, \exists\,K\in\N,C>0''$ holds. 
		\end{itemize}
	\end{theorem}
	The proof of Theorem \ref{thm: main2} will be given in Section \ref{sec: proof}. Interestingly, Theorem \ref{thm: main2} is  somewhat of a different nature than the above mentioned results from \cite{FW04b, MTV90,Wengenroth11}  in the sense that the characterizing properties on the shifted fundamental solutions $E_{\xi,\varepsilon}$ are not only about the behavior of $E_{\xi,\varepsilon}$ on the set $\overline{X}_n$ but also on the larger set $\overline{X}_k$.  In this regard, we mention that $P(D)$ is surjective on $\mathscr{E}(X)$, respectively, $\mathscr{D}'(X)$ if and only if \eqref{eq: Omega}, respectively, \eqref{eq: POmega} without the assumption  $\|E_{\xi,\varepsilon}\|_{\overline{X_k},K}^*\leq \frac{C}{\varepsilon^s}$ holds \cite{Wengenroth11}. Under the a priori assumption that $P(D)$ is surjective on $\mathscr{D}'(X)$, conditions \eqref{eq: Omega} and \eqref{eq: POmega} are equivalent, see \cite[Theorem 1.2]{DeKa22}.
	
	It would be interesting to evaluate the conditions in Theorem \ref{thm: main2} in specific cases in order to obtain concrete necessary and sufficient conditions on $X$ and $P$ for $\mathscr{E}_P(X)$ to satisfy $(\Omega)$ and for $\mathscr{D}'_P(X)$ to satisfy  $(P\Omega)$ and $(P\overline{\overline{\Omega}})$  (cf.\ \cite{MTV90, Langenbruch1998}). As a first result in this direction, we show in Section \ref{sect:convex} that $\mathscr{E}_P(X)$  satisfies $(\Omega)$ for any differential operator $P(D)$ and any open convex set  $X$ by combining Theorem \ref{thm: main2}(a) with  a powerful method to construct fundamental solutions due to H\"ormander \cite[Proof of Theorem 7.3.2]{HoermanderPDO1}. As mentioned above, this result is originally due to Petzsche \cite{Petzsche1980}, who proved it with the aid of the fundamental principle of Ehrenpreis. A completely different proof was recently given by the authors in \cite{DeKa22}. 
	
	Finally, we would like to point out that Theorem \ref{thm: main2} implies that surjectivity of $P(D)$ on $\mathscr{E}(X)$ and $(\Omega)$ for $\mathscr{E}_P(X)$ as well as surjectivity of $P(D)$ on $\mathscr{D}'(X)$ and $(P\Omega)$, respectively, $(P\overline{\overline{\Omega}})$, for $\mathscr{D}'_P(X)$ are preserved under taking finite intersections of open sets:
	
	\begin{corollary}
		Let $P\in\C[\xi_1,\ldots,\xi_d]$ and let $X_1,\ldots,X_n\subseteq\R^d$ be open. Set $X=\bigcap_{j=1}^n X_j$. Suppose $P(D): \mathscr{E}(X_j) \to \mathscr{E}(X_j)$ is surjective and  $\mathscr{E}_P(X_j)$ satisfies $(\Omega)$ for all $1\leq j\leq n$. Then, $P(D): \mathscr{E}(X) \to \mathscr{E}(X)$ is surjective and $\mathscr{E}_P(X)$ satisfies $(\Omega)$ as well. A similar statement holds for $(P\Omega)$ and $(P\overline{\overline{\Omega}})$.
	\end{corollary}
	
	For $(P\Omega)$ the above also follows from \cite[Proposition 8.3]{B-D2006} and the fact that surjectivity of $P(D)$ is preserved under taking finite intersections. However, for $(\Omega)$ and  $(P\overline{\overline{\Omega}})$ we do not see how this may be shown without Theorem \ref{thm: main2}.
	
	\section{Linear topological invariants}\label{sec: preliminaries}
	In this preliminary section we introduce the linear topological invariants $(\Omega)$ for Fr\'echet spaces and $(P\Omega)$ and $(P\overline{\overline{\Omega}})$  for (PLS)-spaces. We refer to \cite{B-D2006,Domanski,M-V} for more information about these conditions and examples of spaces satisfying them.
	
	Throughout, we use standard notation from functional analysis \cite{M-V} and distribution theory \cite{HoermanderPDO1, Schwartz}.
	In particular, given a locally convex space $E$,  we denote by $\mathscr{U}_0(E)$ the filter basis of absolutely convex  neighborhoods of $0$ in $E$ and by $\mathscr{B}(E)$ the family of all absolutely convex bounded sets in $E$.
	
	\subsection{Projective spectra} A \emph{projective spectrum (of locally convex spaces)}  $$\mathcal{E} = (E_n, \varrho_{n+1}^n)_{n \in \N}$$ consists of locally convex spaces $E_n$  and  continuous linear 
	maps $\varrho^n_{n+1}: E_{n+1} \rightarrow E_n$, called the \emph{spectral maps}. We define $\varrho^n_n  = \operatorname{id}_{E_n}$ and $\varrho^n_{m} = \varrho_{n+1}^n \circ \cdots \circ \varrho^{m-1}_{m} : E_m \rightarrow E_n$ for $n,m \in \N$ with $m > n$.  The \emph{projective limit} of $\mathcal{E}$ is defined as
	$$
	\operatorname{Proj} \mathcal{E} = \left\{(x_n)_{n \in \N} \in \prod_{n \in \N} E_n \, | \,  x_n = \varrho^n_{n+1}(x_{n+1}),  \forall n \in \N\right\}.
	$$
	For $n \in \N$ we  write $\varrho^n : \operatorname{Proj} \mathcal{E} \rightarrow E_n, \, (x_j)_{j \in \N} \mapsto x_n$.  We always  endow $\operatorname{Proj} \mathcal{E}$ with its natural projective limit topology. For a projective spectrum $\mathcal{E} = (E_n, \varrho_{n+1}^n)_{n \in \N}$ of Fr\'echet spaces, the projective limit $\operatorname{Proj} \mathcal{E}$ is again a Fr\'echet space. We will implicitly make use of the derived projective limit  $\operatorname{Proj}^1 \mathcal{E}$.  We refer to \cite[Sections 2 and 3]{Wengenroth03} for more information. In particular, see \cite[Theorem 3.1.4]{Wengenroth03} for an explicit definition of  $\operatorname{Proj}^1 \mathcal{E}$.
	
	\subsection{The condition $(\Omega)$ for Fr\'echet spaces}
	A Fr\'echet space $E$ is said to satisfy the condition $(\Omega)$ \cite{M-V} if
	\begin{eqnarray*}
		\begin{gathered}
			\forall U \in \mathscr{U}_0(E) \, \exists V \in \mathscr{U}_0(E) \, \forall W \in \mathscr{U}_0(E) \, \exists s,C > 0 \, \forall \varepsilon \in (0,1) \, : \\
			V \subseteq \varepsilon U + \frac{C}{\varepsilon^s} W.
		\end{gathered}
	\end{eqnarray*} 
	The following result will play a key role in the proof of Theorem \ref{thm: main2}(a).
	
	\begin{lemma}\label{Omega-step} \cite[Lemma 2.4]{DeKa22}
		Let $\mathcal{E} = (E_n, \varrho^n_{n+1})_{n \in \N}$ be a projective spectrum of Fr\'echet spaces. Then,  $\operatorname{Proj}^1 \mathcal{E}  =0$ and  $\operatorname{Proj} \mathcal{E}$ satisfies $(\Omega)$ if and only if
		\begin{eqnarray}\label{Omega-step-form}
			\begin{gathered}
				\forall n \in \N, U \in \mathscr{U}_0(E_n) \, \exists m \geq n, V \in \mathscr{U}_0(E_m) \, \forall k \geq m, W \in \mathscr{U}_0(E_k) \\ 
				\exists s,C > 0 \, \forall \varepsilon \in (0,1) \, : 
				\varrho^n_m(V) \subseteq \varepsilon U + \frac{C}{\varepsilon^s} \varrho^n_k(W). 
			\end{gathered}
		\end{eqnarray}
	\end{lemma}
	
	\subsection{The conditions $(P\Omega)$ and $(P\overline{\overline{\Omega}})$  for (PLS)-spaces} A locally convex space $E$ is called a (PLS)-space if it can be written as the projective limit of a spectrum of (DFS)-spaces. 
	
	Let  $\mathcal{E} = (E_n, \varrho^n_{n+1})_{n \in \N}$ be a spectrum of (DFS)-spaces. We call $\mathcal{E}$ \emph{strongly reduced} if
	\[
	\forall n \in \N \, \exists m \geq n \, : \, \varrho^n_m(E_m) \subseteq \overline{\varrho^n( \operatorname{Proj} \mathcal{E} )}.
	\] 
	The spectrum $\mathcal{E}$ is said to satisfy $(P\Omega)$ if
	\begin{eqnarray}\label{Pomega-abstract}
		\begin{gathered}
			\forall n \in \N \, \exists m \geq n \, \forall k \geq m \, \exists B \in \mathscr{B}(E_n) \, \forall M \in \mathscr{B}(E_m) \, \exists K \in \mathscr{B}(E_k),s,C > 0 \\
			\forall \varepsilon \in (0,1) \, :
			\varrho^n_m(M) \subseteq \varepsilon B + \frac{C}{\varepsilon^s} \varrho^n_k(K).
		\end{gathered}
	\end{eqnarray}
	The spectrum $\mathcal{E}$ is said to satisfy $(P\overline{\overline{\Omega}})$ if  \eqref{Pomega-abstract} with ``$ \exists K \in \mathscr{B}(E_k), s,C > 0$'' replaced by ``$\forall s > 0 \, \exists K \in \mathscr{B}(E_k), C > 0$'' holds.
	
	A (PLS)-space $E$ is said to satisfy $(P\Omega)$, respectively, $(P\overline{\overline{\Omega}})$  if  $E = \operatorname{Proj}\mathcal{E}$ for some strongly reduced spectrum $\mathcal{E}$ of (DFS)-spaces that satisfies $(P\Omega)$, respectively, $(P\overline{\overline{\Omega}})$. This notion is well-defined as \cite[Proposition 3.3.8]{Wengenroth03} yields that all strongly reduced projective spectra $\mathcal{E}$ of (DFS)-spaces with $E = \operatorname{Proj}\mathcal{E}$ are equivalent (in the sense of  \cite[Definition 3.1.6]{Wengenroth03}). The bipolar theorem and \cite[Lemma 4.5]{B-D2006} imply that the above definitions of $(P\Omega)$ and $(P\overline{\overline{\Omega}})$ are equivalent to the original ones from \cite{B-D2006}.
	
	\section{Proof of Theorem \ref{thm: main2}}\label{sec: proof}
	This section is devoted to the proof of Theorem \ref{thm: main2}.  We fix $P\in\C[\xi_1,\ldots,\xi_d]\backslash\{0\}$, an open set $X\subseteq\R^d$, and an exhaustion  by relatively compact open subsets $(X_n)_{n\in\N}$ of $X$. For $n, N \in \N$ we write  $\| \, \cdot \, \|_{n,N} =  \| \, \cdot \, \|_{\overline{X_n},N}$ and $\| \, \cdot \, \|^*_{n,N} =  \| \, \cdot \, \|^*_{\overline{X_n},N}$.
	For $\xi\in\R^d$ and $r>0$ we denote by $B(\xi,r)$ the open ball in $\R^d$ with center $\xi$ and radius $r$. Moreover, for $p \in\{1,\infty\}$ and $N \in \N$ we set
	$$
	\| \varphi\|_{L^p,N} = \max_{|\alpha| \leq N} \| \varphi^{(\alpha)}\|_{L^p}, \qquad \varphi \in \mathscr{D}(\R^d).
	$$
	We fix $\chi\in\mathscr{D}(B(0,1))$ with $\chi\geq 0$ and $\int_{\R^d}\chi(x)\dx=1$, and set  $\chi_\varepsilon(x)=\varepsilon^{-d}\chi(x/\varepsilon)$ for $\varepsilon > 0$.

	\subsection{Proof of Theorem \ref{thm: main2}(a)}\label{sec: proof1} We write $\mathscr{E}(X)$ for the space of smooth functions in $X$ endowed with its natural Fr\'echet space topology. We set
	$$
	\mathscr{E}_P(X) = \{ f \in \mathscr{E}(X) \, | \, P(D)f = 0 \}
	$$
	and endow it with the relative topology induced by $\mathscr{E}(X)$.

	Let $n \in \N$. We write $\mathscr{E}(\overline{X_n})$ for the space of continuous functions in $\overline{X_n}$ whose restrictions to $X_n$ are smooth such that all partial derivatives extend continuously to $\overline{X_n}$, endowed with its natural Fr\'echet space topology, i.e, the one induced by the sequence of norms  $( \| \, \cdot \, \|_{n,N})_{N \in \N}$. We define
	$$
	\mathscr{E}_P(\overline{X_n}) = \{ f \in \mathscr{E}(\overline{X_n}) \, | \, P(D)f = 0 \}
	$$
	and endow it with the relative topology induced by $\mathscr{E}(\overline{X_n})$. Since $\mathscr{E}_P(\overline{X_n})$ is closed in $\mathscr{E}(\overline{X_n})$,  it is a Fr\'echet space. For $N \in \N$ we set
	$$
	U_{n,N} = \{ f \in \mathscr{E}_P(\overline{X_n}) \, | \, \| f \|_{n,N} \leq 1 \}.
	$$
	Note that $\left(\frac{1}{N+1}U_{n,N}\right)_{N \in \N}$ is a decreasing fundamental sequence of absolutely convex neighborhoods of $0$ in $\mathscr{E}(\overline{X_n})$. 
	
	Consider the projective spectrum $(\mathscr{E}_P(\overline{X_n}),\varrho_{n+1}^n )_{n \in \N}$ with  $\varrho_{n+1}^n$ the restriction map from $\mathscr{E}_P(\overline{X_{n+1}})$ to $\mathscr{E}_P(\overline{X_n})$. Then,
	$$
	\mathscr{E}_P (X) =\operatorname{Proj}(\mathscr{E}_P(\overline{X_n}),\varrho_{n+1}^n )_{n \in \N}.
	$$ 
	By \cite[Lemma 3.1$(i)$]{DeKa22} (see also \cite[Section 3.4.4]{Wengenroth03}), $P(D): \mathscr{E}(X) \to \mathscr{E}(X)$ is surjective if and only if 
	$$
	\operatorname{Proj}^1(\mathscr{E}_P(\overline{X_n}),\varrho_{n+1}^n )_{n \in \N} = 0.
	$$
	Hence, Lemma \ref{Omega-step} and a simple rescaling argument yield that  $P(D): \mathscr{E}(X) \to \mathscr{E}(X)$ is surjective and  $\mathscr{E}_P(X)$ satisfies $(\Omega)$ if and only if 
	\begin{eqnarray}\label{Omega-PDO}
		\begin{gathered}
			\forall n,N \in \N \, \exists m \geq n, M \geq N  \, \exists k \geq m, K \geq M \, \exists s,C> 0 \, \forall \varepsilon \in (0,1)\, : \\ 
			U_{m,M}\subseteq \varepsilon U_{n,N} + \frac{C}{\varepsilon^s} U_{k,K},
		\end{gathered}
	\end{eqnarray}
	where we did not write the restriction maps explicitly, as we shall not do in the sequel either. We are ready to show Theorem \ref{thm: main2}(a).
	\subsubsection*{Sufficiency of  \eqref{eq: Omega}} It suffices to show \eqref{Omega-PDO}. Let $n,N \in \N$ be arbitrary. Choose $\widetilde{m}, \widetilde{N}$ according to \eqref{eq: Omega} for $n+1$. Set  $m=\widetilde{m}+1$ and $M=N+\widetilde{N}+\deg P+1$. Let $k \geq m$, $K \geq M$ be arbitrary. Choose $\psi \in \mathscr{D}(X_m)$ such that  $\psi=1$ in a neighborhood of $\overline{X_{\widetilde{m}}}$. Pick $\varepsilon_0 \in (0,1]$ such that
	$\psi = 1$ on $\overline{X_{\widetilde{m}}} + B(0,\varepsilon_0)$, $\supp \psi + B(0,\varepsilon_0) \subseteq X_{m}$, $\overline{X_n} + B(0,\varepsilon_0) \subseteq X_{n+1}$, $\overline{X_k} + B(0,\varepsilon_0) \subseteq X_{k+1}$. 
	Cover the compact set $\overline{X_k}\backslash X_{\widetilde{m}}$ by finitely many balls $B(\xi_j,\varepsilon_0)$, $j\in J$, with $\xi_j \in X \backslash \overline{X_{\widetilde{m}}}$, and choose $\varphi_j \in \mathscr{D}(B(\xi_j,\varepsilon_0))$, $j \in J$, such that 
	$\sum_{j \in J} \varphi_j = 1$ in a neighborhood of $\overline{X_k}\backslash X_{\widetilde{m}}$. As $J$ is finite, \eqref{eq: Omega} for $k+1$ implies that there are $\widetilde{K}\in\N, \widetilde{s},\widetilde{C}>0$ such that for all $\varepsilon \in (0,\varepsilon_0)$ there exist $E_{\xi_j,\varepsilon}\in\mathscr{D}'(\R^d)$, $j \in J$, with  $P(D)E_{\xi_j,\varepsilon}=\delta_{\xi_j}$ in $X_{k+1}$  such that  
	\begin{equation}
		\label{shiffted-ineq}
		\|E_{\xi_j,\varepsilon}\|_{n+1,\tilde{N}}^*\leq \varepsilon \quad \mbox{and} \quad \|E_{\xi_j,\varepsilon}\|_{k+1,\widetilde{K}}^*\leq \frac{\widetilde{C}}{\varepsilon^{\widetilde{s}}}. 
	\end{equation}  
	Let $f\in\mathscr{E}_P(\overline{X_m})$ be arbitrary. For $\varepsilon \in (0,\varepsilon_0)$ we define $f_\varepsilon = (\psi f)\ast \chi_\varepsilon \in \mathscr{E}_P(\overline{X_{\widetilde{m}}})$ and
	$$
	h_\varepsilon=\sum_{j\in J}E_{\xi_j,\varepsilon}\ast\delta_{-\xi_j}\ast (\varphi_jP(D) f_\varepsilon).
	$$
	Since $\delta_{-\xi_j}\ast  (\varphi_jP(D) f_\varepsilon)= (\varphi_jP(D) f_\varepsilon)(\cdot + \xi_j)\in\mathscr{D}(B(0,\varepsilon_0))$, $j \in J$, it holds that $P(D)h_\varepsilon=\sum_{j\in J} \varphi_jP(D) f_\varepsilon$ in a neighborhood of $\overline{X_k}$. As $\sum_{j\in J} \varphi_j = 1$ in a neighborhood of $\overline{X_k}\backslash X_{\widetilde{m}}$ and $f_\varepsilon  \in \mathscr{E}_P(\overline{X_{\widetilde{m}}})$, we obtain that $P(D)h_\varepsilon=P(D) f_\varepsilon$ in a neighborhood of $\overline{X_k}$ and thus  $h_\varepsilon  \in \mathscr{E}_P(\overline{X_{\widetilde{m}}})$ and  $f_\varepsilon-h_\varepsilon\in\mathscr{E}_P(\overline{X_k})$. We decompose $f$ as follows
	$$
	f = (f - f_\varepsilon + h_\varepsilon) + (f_\varepsilon - h_\varepsilon) \in \mathscr{E}_P(\overline{X_n}) + \mathscr{E}_P(\overline{X_k}).
	$$
	We claim that there are $s,C_i >0$, $i = 1,2,3,4$, such that for all  $f\in\mathscr{E}_P(\overline{X_m})$ and $\varepsilon \in (0,\varepsilon_0)$
	\begin{eqnarray*}
		\begin{gathered}
			\|f-f_\varepsilon \|_{n,N} \leq C_1\varepsilon \| f \|_{m,M}, \qquad \|h_\varepsilon \|_{n,N} \leq C_2\varepsilon \| f \|_{m,M}, \\
			\|f_\varepsilon \|_{k,K} \leq  \frac{C_3}{\varepsilon^s} \| f \|_{m,M}, \qquad \|h_\varepsilon \|_{k,K} \leq  \frac{C_4}{\varepsilon^s} \| f \|_{m,M},
		\end{gathered}
	\end{eqnarray*}
	which implies  \eqref{Omega-PDO}. Let $f\in\mathscr{E}_P(\overline{X_m})$ and $\varepsilon \in (0,\varepsilon_0)$ be arbitrary. By the mean value theorem, we find that 
	$$
	\|f- f_{\varepsilon}\|_{n,N} \leq \varepsilon \sqrt{d}  \| f \|_{n+1,N+1} \leq \varepsilon \sqrt{d}  \| f \|_{m,M}.
	$$
	Furthermore, it holds that
	$$
	\|f_\varepsilon \|_{k,K} \leq  \frac{\| \chi\|_{L^1,K} } {\varepsilon^{K}} \|\psi f \|_{L^\infty} \leq \frac{\| \chi\|_{L^1,K} \| \psi\|_{L^\infty} } {\varepsilon^{K}} \| f \|_{m,M}.
	$$ 
	By the first inequality in \eqref{shiffted-ineq}, we obtain that
	\begin{eqnarray*}
		\|h_\varepsilon\|_{n,N}&\leq& \sum_{j\in J} \|E_{\xi_j,\varepsilon}\|_{n+1, \widetilde{N}}^*\|(\varphi_jP(D) f_\varepsilon)(\cdot + \xi_j)\|_{L^\infty, N+\widetilde{N}}\\
		&\leq &\varepsilon\sum_{j\in J}\|\varphi_j ((P(D)(\psi f))\ast \chi_\varepsilon)\|_{L^\infty,N+\widetilde{N}}\\
		&\leq & C'_2 \varepsilon \|P(D)(\psi f)\|_{L^\infty, N+\widetilde{N}}\\
		&\leq & C_2 \varepsilon \|f\|_{m,N+\widetilde{N} + \deg P}\leq  C_2 \varepsilon \|f\|_{m,M},
	\end{eqnarray*}
	for some $C'_2,C_2 > 0$.
	Similarly, by the second inequality in \eqref{shiffted-ineq}, we find that
	\begin{eqnarray*}
		\|h_\varepsilon\|_{k,K}&\leq& \sum_{j\in J} \|E_{\xi_j,\varepsilon}\|_{k+1, \widetilde{K}}^*\|(\varphi_jP(D) f_\varepsilon)(\cdot + \xi_j)\|_{L^\infty,K+\widetilde{K}}\\
		&\leq &\frac{\widetilde{C}}{\varepsilon^{\widetilde{s}}} \sum_{j\in J}\|\varphi_j ((P(D)(\psi f))\ast \chi_\varepsilon)\|_{L^\infty,K+\widetilde{K}}\\
		&\leq & \frac{C'_4}{\varepsilon^{\widetilde{s}}} \|P(D)(\psi f)\|_{L^\infty} \| \chi_\varepsilon \|_{L^1,K+\widetilde{K}}\\
		&\leq & \frac{C_4}{\varepsilon^{\widetilde{s} + K + \widetilde{K}}} \|f\|_{m, \deg P} \leq   \frac{C_4}{\varepsilon^{\widetilde{s} + K + \widetilde{K}}} \|f\|_{m, M}, 
	\end{eqnarray*}
	for some $C'_4,C_4 > 0$.This proves the claim with $s = \widetilde{s} + K + \widetilde{K}$. \hfill $\qed$

	\subsubsection*{Necessity of  \eqref{eq: Omega}} As explained above, condition \eqref{Omega-PDO} holds. Let $F \in \mathscr{D}'(\R^d)$ be a fundamental solution for $P(D)$ of finite order $q$. Let $n \in \N$ be arbitrary. Choose $m, \widetilde{M} \in \N$ according to \eqref{Omega-PDO} for $n$ and  $0$. Set $N = q+1$. Let $k \geq m$ and $\xi \in X \backslash \overline{X_m}$ be arbitrary. Set $K = q+1$. \eqref{Omega-PDO} for $k+1$ and $0$ implies that there are $\widetilde{C}, \widetilde{s} > 0$ such that for all $\delta \in (0,1)$ 
	\begin{equation}
		\label{omega-delta}
		U_{m,\widetilde{M}}\subseteq \delta U_{n,0} + \frac{\widetilde{C}}{\delta^{\widetilde{s}}} U_{k+1,0}.
	\end{equation}
	Let $\varepsilon_0 \in (0,1]$ be such that $B(\xi, \varepsilon_0) \subseteq X \backslash \overline{X_m}$. 
	Set $F_\xi = F \ast \delta_\xi \in \mathscr{D}'(\R^d)$. For all $\varepsilon \in (0,\varepsilon_0)$ it holds that $F_\xi\ast \chi_\varepsilon\in\mathscr{E}_P(\overline{X_{m}})$ and
	$$
	\| F_\xi\ast \chi_\varepsilon \|_{m,\widetilde{M}} \leq \frac{C'}{\varepsilon^{d+\widetilde{M}+q}}
	$$ 
	with $C' = \| F_\xi \|^*_{\overline{X_m +B(0,\varepsilon_0)},q} $. Hence, \eqref{omega-delta} with $\delta = \varepsilon^{d+\widetilde{M}+q  + 1}$ implies that
	$$
	F_\xi\ast \chi_\varepsilon\in  \frac{C'}{\varepsilon^{d+\widetilde{M}+q}}U_{m,\widetilde{M}}\subseteq C' \varepsilon U_{n,0} + \frac{C'\widetilde{C}}{\varepsilon^{s}} U_{k+1,0},
	$$
	with $s = d+\widetilde{M}+q + \widetilde{s}(d+\widetilde{M}+q  + 1)$. Let $f_{\xi,\varepsilon}\in C' \varepsilon U_{n,0}$ and $h_{\xi,\varepsilon}\in C'\widetilde{C}\varepsilon^{-s} U_{k+1,0}$ be such that
	\begin{equation}\label{eq: decomposition}
		F_\xi\ast\chi_\varepsilon=f_{\xi,\varepsilon}+h_{\xi,\varepsilon}.
	\end{equation}
	Choose  $\psi \in\mathscr{D}(X_{k+1})$ such that $\psi =1$ in a neighborhood of $\overline{X_k}$ and define $E_{\xi,\varepsilon}=F_\xi-\psi h_{\xi,\varepsilon}\in\mathscr{D}'(\R^d)$. Then,  $P(D)E_{\xi,\varepsilon}=\delta_\xi$ in $X_k$. Moreover, for all $\varepsilon \in (0,\varepsilon_0)$ it holds that 
	\begin{eqnarray*}
		\|E_{\xi,\varepsilon}\|_{n,q+1}^\ast&\leq&\|F_{\xi}-F_\xi\ast\chi_\varepsilon\|_{n,q+1}^\ast+\|F_\xi\ast\chi_\varepsilon -\psi h_{\xi,\varepsilon}\|_{n,q+1}^\ast\\
		&\leq &\| F_{\xi} \|^*_{\overline{X_n +B(0,\varepsilon_0)},q} \sqrt{d} \varepsilon +\|f_{\xi,\varepsilon}\|_{n,0}\\
		&\leq &(\| F_{\xi} \|^*_{\overline{X_n +B(0,\varepsilon_0)},q} \sqrt{d}+C')\varepsilon,
	\end{eqnarray*}
	where we used the mean value theorem, and
	\begin{eqnarray*}
		\|E_{\xi,\varepsilon}\|_{k,q+1}^\ast&\leq&\| F_\xi\|_{k,q+1}^\ast + \| \psi h_{\xi,\varepsilon}\|_{k,q+1}^\ast \\
		&\leq& \| F_\xi\|_{k,q+1}^\ast + |X_k|\|h_{\xi,\varepsilon}\|_{k,0} \\
		&\leq& \frac{\| F_\xi\|_{k,q+1}^\ast + C'\widetilde{C}|X_k|}{\varepsilon^s},
	\end{eqnarray*}
	where $|X_k|$ denotes the Lebesgue measure of $X_k$. This completes the proof. \hfill \qed
	
	\subsection{Proof of Theorem \ref{thm: main2}(b) and (c)}\label{sec: proof2} We write $\mathscr{D}'(X)$ for the space of distributions in $X$ endowed with its strong dual topology. We set
	$$
	\mathscr{D}'_P(X) = \{ f \in \mathscr{D}'(X) \, | \, P(D)f = 0 \}
	$$
	and endow it with the relative topology induced by $\mathscr{D}'(X)$.
	
	In \cite[Theorem $(5)$]{Wengenroth11} it is shown that the mapping $P(D): \mathscr{D}'(X) \to \mathscr{D}'(X)$ is surjective if and only if 
	\begin{eqnarray}\label{surj-dist}
		\begin{gathered}
			\forall\,n\in\N\,\exists\,m\geq n\,\forall\,k\geq m, N\in\N, \xi\in X\backslash \overline{X_m},\varepsilon \in (0,1)\\
			\exists\,E_{\xi,\varepsilon}\in\mathscr{D}'(\R^d)\cap C^N(\overline{X_n}) \mbox{ with }  P(D)E_{\xi,\varepsilon}=\delta_\xi\mbox{ in }X_k  \mbox{ such that }   \\
			\|E_{\xi,\varepsilon}\|_{\overline{X_n},N}\leq \varepsilon.
		\end{gathered}
	\end{eqnarray}
	Let $n \in \N$. We endow the space $\mathscr{D}_{\overline{X_n}}$ of smooth functions with support in  $\overline{X_n}$ with the relative topology induced by $\mathscr{E}(\overline{X_n})$. We write $\mathscr{D}'(\overline{X_n})$ for the strong dual of $\mathscr{D}_{\overline{X_n}}$. Then, $\mathscr{D}'(\overline{X_n})$ is a (DFS)-space. We define
	$$
	\mathscr{D}'_P(\overline{X_n}) = \{ f \in \mathscr{D}'(\overline{X_n}) \, | \, P(D)f = 0 \}
	$$
	and endow it with the relative topology induced by $\mathscr{D}'(\overline{X_n})$. Since $\mathscr{D}'_P(\overline{X_n})$ is closed in $\mathscr{D}'(\overline{X_n})$,  it is a  (DFS)-space. For $N \in \N$ we set
	$$
	B_{n,N} = \{ f \in \mathscr{D}'_P(\overline{X_n}) \, | \, \| f \|^*_{n,N} \leq 1 \}.
	$$
	Note that $(NB_{n,N})_{N \in \N}$ is an increasing fundamental sequence of absolutely convex bounded sets  in $\mathscr{D}_P'(\overline{X_n})$. 
	
	Consider the projective spectrum $(\mathscr{D}'_P(\overline{X_n}),\varrho_{n+1}^n )_{n \in \N}$ with  $\varrho_{n+1}^n$ the restriction map from $\mathscr{D}'_P(\overline{X_{n+1}})$ to $\mathscr{D}'_P(\overline{X_n})$. Then,
	$$
	\mathscr{D}'_P (X) =\operatorname{Proj}(\mathscr{D}'_P(\overline{X_n}),\varrho_{n+1}^n )_{n \in \N}.
	$$ 
	By \cite[Lemma 3.1$(ii)$]{DeKa22} (see also \cite[Section 3.4.5]{Wengenroth03}), $P(D): \mathscr{D}'(X) \to \mathscr{D}'(X)$ is surjective if and only if 
	$$
	\operatorname{Proj}^1(\mathscr{D}'_P(\overline{X_n}),\varrho_{n+1}^n )_{n \in \N} = 0.
	$$
	The latter condition implies that $(\mathscr{D}'_P(\overline{X_n}),\varrho_{n+1}^n )_{n \in \N}$ is strongly reduced \cite[Theorem 3.2.9]{Wengenroth03}.  Hence, if  $P(D): \mathscr{D}'(X) \to \mathscr{D}'(X)$ is surjective, $\mathscr{D}'_P (X)$ satisfies $(P\Omega)$, respectively,  $(P\overline{\overline{\Omega}})$ if and only if $(\mathscr{D}'_P(\overline{X_n}),\varrho_{n+1}^n )_{n \in \N}$ does so. A simple rescaling argument yields that  $(\mathscr{D}'_P(\overline{X_n}),\varrho_{n+1}^n )_{n \in \N}$ satisfies $(P\Omega)$ if and only if
	\begin{eqnarray}\label{POmega-PDO}
		\begin{gathered}
			\forall n \in \N \, \exists m \geq n \, \forall k \geq m \, \exists N \in \N \, \forall M \in \N \, \exists K \in \N, s,C > 0 \, \forall \varepsilon \in (0,1) \, : \\ 
			B_{m,M} \subseteq \varepsilon B_{n,N} + \frac{C}{\varepsilon^s} B_{k,K},
		\end{gathered}
	\end{eqnarray}
	where, as before, we did not write the restriction maps explicitly.  Similarly, the spectrum $(\mathscr{D}'_P(\overline{X_n}),\varrho_{n+1}^n )_{n \in \N}$ satisfies $(P\overline{\overline{\Omega}})$  if and only if   \eqref{POmega-PDO} with ``$ \exists K \in \N, s,C > 0$'' replaced by ``$\forall s > 0 \, \exists K \in \N, C > 0$'' holds.  We now show Theorem \ref{thm: main2}(b).
	
	\subsubsection*{Sufficiency of  \eqref{eq: POmega}}  Clearly,  \eqref{eq: POmega}  implies \eqref{surj-dist} and thus that $P(D): \mathscr{D}'(X) \to \mathscr{D}'(X)$ is surjective. Hence, by the above discussion, it suffices to show \eqref{POmega-PDO}. Let $n \in \N$ be arbitrary. Choose $m$ according to \eqref{eq: POmega} for $n+1$. Let $k \geq m$ be arbitrary. Set $N = 0$. Let $M \in \N$ be arbitrary. Pick $\varepsilon_0 \in (0,1]$ such that $\overline{X_n} + B(0,\varepsilon_0) \subseteq X_{n+1}$ and $\overline{X_k} + B(0,\varepsilon_0) \subseteq X_{k+1}$.  Cover the compact set $\overline{X_k}\backslash X_{m}$ by finitely many balls $B(\xi_j,\varepsilon_0)$, $j\in J$, with $\xi_j \in X \backslash \overline{X_m}$, and choose $\varphi_j \in \mathscr{D}(B(\xi_j,\varepsilon_0))$, $j \in J$, such that 
	$\sum_{j \in J} \varphi_j = 1$ in a neighborhood of $\overline{X_k}\backslash X_{m}$. As $J$ is finite, \eqref{eq: POmega} for $k+1$ and $M + \deg P$  implies that there are $\widetilde{K}\in\N, \widetilde{s},\widetilde{C}>0$ such that for all $\varepsilon \in (0,\varepsilon_0)$ there exist $E_{\xi_j,\varepsilon}\in\mathscr{D}'(\R^d) \cap C^{M + \deg P}(\overline{X_{n+1}})$, $j \in J$, with  $P(D)E_{\xi_j,\varepsilon}=\delta_{\xi_j}$ in $X_{k+1}$  such that  
	\begin{equation}
		\label{shiffted-ineq-2}
		\|E_{\xi_j,\varepsilon}\|_{\overline{X_{n+1}},M+\deg P}\leq \varepsilon \quad \mbox{and} \quad \|E_{\xi_j,\varepsilon}\|_{k+1,\widetilde{K}}^*\leq \frac{\widetilde{C}}{\varepsilon^{\widetilde{s}}}. 
	\end{equation}  
	Pick $\psi\in\mathscr{D}(X_{m})$ with $\psi=1$ in a neighborhood of $\overline{X_n}$. Let $f\in \mathscr{D}'_P(\overline{X_m})$ with $\|f\|^*_{m,M} < \infty$ be arbitrary. For $\varepsilon \in (0,\varepsilon_0)$ we define 
	$$
	h_\varepsilon=\sum_{j\in J}E_{\xi_j,\varepsilon}\ast\delta_{-\xi_j}\ast (\varphi_jP(D)(\psi f)).
	$$
	By the same reasoning as in the proof of part (a) it follows that $h_\varepsilon \in \mathscr{D}'_P(\overline{X_n})$ and $\psi f - h_\varepsilon \in \mathscr{D}'_P(\overline{X_k})$. Furthermore, as $E_{\xi_j,\varepsilon}\in\mathscr{D}'(\R^d) \cap C^{M + \deg P}(\overline{X_{n+1}})$ and the distributions $\delta_{-\xi_j}\ast (\varphi_jP(D)(\psi f)) =  \varphi_jP(D) (\psi f)(\cdot + \xi_j)$, $j \in J$, have order at most $M + \deg P$ and are supported in $B(0,\varepsilon_0)$, it holds that $h_\varepsilon \in C(\overline{X_n})$. We decompose $f$ as follows  in $X_n$
	$$
	f = \psi f =  h_\varepsilon + (\psi f - h_\varepsilon) \in  ( \mathscr{D}'_P(\overline{X_n}) \cap C(\overline{X_n})  ) +  \mathscr{D}'_P(\overline{X_k}).
	$$
	We claim that there are $K \in \N$, $s,C_1,C_2 >0$ such that for all  $f\in \mathscr{D}'_P(\overline{X_m})$ with $\|f\|^*_{m,M} < \infty$ and $\varepsilon \in (0,\varepsilon_0)$
	$$
	\| h_\varepsilon \|^*_{n,0} \leq  C_1\varepsilon \| f \|^*_{m,M}, \qquad  \|\psi f - h_\varepsilon \|^*_{k,K} \leq  \frac{C_2}{\varepsilon^s} \| f \|^*_{m,M},
	$$
	which implies  \eqref{POmega-PDO}. Let $f\in \mathscr{D}'_P(\overline{X_m})$ with $\|f\|^*_{m,M} < \infty$ and $\varepsilon \in (0,\varepsilon_0)$ be arbitrary. Choose $\rho \in\mathscr{D}(X_{m})$ with $\rho=1$ in a neighborhood of  $\supp \psi$. The first inequality in \eqref{shiffted-ineq-2} implies that
	\begin{eqnarray*}
		\| h_\varepsilon \|^*_{n,0}&\leq&  |X_n|\| h_\varepsilon \|_{n,0} \\
		&\leq&  |X_n|\sum_{j\in J} \| P(D)(\psi f)\|^*_{m, M + \deg P} \sup_{x \in \overline{X_n}} \| (\varphi_j\rho) E_{\xi_j,\varepsilon}(x+\xi_j- \cdot) \|_{L^\infty, M + \deg P} \\
		&\leq& C_1\| f\|^*_{m, M} \| E_{\xi_j,\varepsilon}\|_{n+1, M + \deg P} \\
		&\leq & C_1\varepsilon \| f \|^*_{m,M}
	\end{eqnarray*}
	for some $C_1 > 0$. Next, by the second inequality in \eqref{shiffted-ineq-2}, we obtain that   for all $\varphi\in\mathscr{D}_{\overline{X_k}}$ 
	\begin{eqnarray*}
		|\langle h_\varepsilon, \varphi \rangle| &\leq& \sum_{j \in J}|\langle E_{\xi_j,\varepsilon}\ast\delta_{-\xi_j}\ast (\varphi_j P(D)(\psi f)),\varphi\rangle| \\
		&=&\sum_{j \in J}|\langle E_{\xi_j,\varepsilon},(\delta_{-\xi_j}\ast (\varphi_j P(D)(\psi f)))^{\vee}\ast\varphi\rangle|\\
		&\leq& \sum_{j \in J}\|E_{\xi_j,\varepsilon}\|_{k+1,\widetilde{K}}^\ast\| (\delta_{-\xi_j}\ast (\varphi_j P(D)(\psi f)))^{\vee}\ast\varphi \|_{L^\infty,\widetilde{K}}\\
		&\leq& \frac{\widetilde{C}}{\varepsilon^{\widetilde{s}}}\sum_{j \in J}   \| P(D)(\psi f)\|^*_{m, M + \deg P} \times\\
		&&\times\sup_{x \in \R^d} \|  (\varphi_j\rho) \varphi(\cdot - \xi_j-x)  \|_{L^\infty,\widetilde{K}+M + \deg P}\\
		&\leq& \frac{C'_2}{\varepsilon^{\widetilde{s}}} \|f\|^*_{m, M} \| \varphi \|_{L^\infty,\widetilde{K}+M + \deg P},
	\end{eqnarray*}
	for some $C'_2 > 0$, whence 
	$$
	\|h_\varepsilon \|^*_{k,\widetilde{K} + M + \deg P} \leq \frac{C'_2}{\varepsilon^{\widetilde{s}}} \|f\|^*_{m, M}.
	$$
	Furthermore,
	$$
	\|\psi f  \|^*_{k,\widetilde{K} + M + \deg P} \leq C''_2\|f\|^*_{m, M},
	$$
	for some $C''_2 > 0$. Therefore,
	$$
	\|\psi f - h_\varepsilon \|^*_{k,\widetilde{K} + M + \deg P}  \leq \frac{C''_2 + C'_2}{\varepsilon^{\widetilde{s}}} \|f\|^*_{m, M}.
	$$
	This shows the claim with $K = \widetilde{K} + M + \deg P$ and $s =\widetilde{s}$. \hfill \qed
	
	\subsubsection*{Necessity of  \eqref{eq: POmega}} As explained above, conditions \eqref{surj-dist} and \eqref{POmega-PDO}  hold.  Let $n \in \N$  be arbitrary. Choose $\widetilde{m}$ according to \eqref{POmega-PDO} for $n+1$ and  $m$ according to \eqref{surj-dist} for $\widetilde{m}$. Let $k \geq m$, $N \in \N$, and $\xi \in X \backslash \overline{X_m}$ be arbitrary.   \eqref{surj-dist}  for $k$ and $N+1$ implies that there exist $F_\xi \in \mathscr{D}'(\R^d)\cap C^{N+1}(\overline{X_{\widetilde{m}}})$ with $P(D)F_\xi=\delta_\xi$ in $X_k$ such that $\|F_\xi\|_{\widetilde{m},N+1}\leq \min\{1, 1/|X_{\widetilde{m}}|\}$.  By \eqref{POmega-PDO} for $k+1$, we obtain that there are $\widetilde{N}, \widetilde{K} \in \N$, $\widetilde{s}, \widetilde{C} >0$ such that for all $\delta \in (0,1)$
	\begin{equation}
		\label{Pomega-inproof}
		B_{\widetilde{m},0} \subseteq \delta B_{n+1,\widetilde{N}} + \frac{\widetilde{C}}{\delta^{\widetilde{s}}} B_{k+1,\widetilde{K}}.
	\end{equation}
	Note that $F_{\xi} \in \mathscr{D}'_P(\overline{X_{\widetilde{m}}})$ and $\|F_\xi\|^*_{\widetilde{m},0} \leq |X_{\widetilde{m}}|\|F_\xi\|_{\widetilde{m},0} \leq 1$, whence $F_{\xi}\in B_{\widetilde{m},0}$. Therefore, \eqref{Pomega-inproof} yields that for all $\delta \in (0,1)$ there are  $G_{\xi,\delta} \in \delta B_{n+1,\widetilde{N}}$ and $H_{\xi,\delta} \in \widetilde{C}\delta^{-\widetilde{s}} B_{k+1,\widetilde{K}}$ such that $F_{\xi} = G_{\xi,\delta} + H_{\xi,\delta}$ in $X_{n+1}$. Let $\psi \in \mathscr{D}(X_{k+1})$ be such that $\psi = 1$ on a neighborhood of $ \overline{X_k}$. Choose $\varepsilon_0 \in (0,1]$ such that $\psi =1$ on $\overline{X_{k}} + B(0,\varepsilon_0)$ and  $\overline{X_n} + B(0,\varepsilon_0) \subseteq X_{n+1}$. For $\delta \in (0,1)$ and $\varepsilon \in (0,\varepsilon_0)$ we define
	$$
	E_{\xi,\varepsilon,\delta} = F_\xi - (\psi H_{\xi,\delta}) \ast \chi_\varepsilon \in \mathscr{D}'(\R^d).
	$$
	Since $H_{\xi,\delta} \in \mathscr{D}'_P(\overline{X_{k+1}})$, we have that $(\psi H_{\xi,\delta}) \ast \chi_\varepsilon \in \mathscr{E}_P(\overline{X_{k}})$. This implies that  $ E_{\xi,\varepsilon,\delta} \in C^{N}(\overline{X_n})$ and $P(D) E_{\xi,\varepsilon,\delta} = \delta_\xi$ in $X_k$. As $\psi = 1$ on $X_{n+1}$, it holds that $(\psi H_{\xi,\delta}) \ast \chi_\varepsilon = H_{\xi,\delta} \ast \chi_\varepsilon$ on $\overline{X_n}$. Hence, we obtain that 
	\begin{eqnarray*}
		\|E_{\xi,\varepsilon,\delta}\|_{n,N}&\leq&\|F_\xi-F_\xi\ast\chi_\varepsilon\|_{n,N}+\|F_\xi\ast\chi_\varepsilon-H_{\xi,\delta}\ast\chi_\varepsilon\|_{n,N}\\\
		&\leq & \sqrt{d}\varepsilon+\|G_{\xi,\delta}\ast\chi_\varepsilon\|_{n,N}\\
		&\leq &\sqrt{d} \varepsilon +\|\chi\|_{L^\infty,N+ \widetilde{N}}\frac{\delta}{\varepsilon^{N+\widetilde{N}+d}},
	\end{eqnarray*}
	where we used the mean value theorem.  Let $r \in \N$ be the order of $F_\xi$ in $\overline{X_k}$ and set $K = \max\{r, \widetilde{K} \}$. Then,
	$$
	\|E_{\xi,\varepsilon,\delta}\|_{k,K}^*\leq \|F_\xi\|_{k,r}^*+\|(\psi H_{\xi,\delta})\ast\chi_\varepsilon\|_{k,\widetilde{K}}^\ast\leq \frac{\|F_\xi\|_{k,r}^*+C'}{\delta^{\widetilde{s}}},
	$$
	for some $C' > 0$. For $\varepsilon \in (0,\varepsilon_0)$ we set $\delta_\varepsilon = \varepsilon^{N+\widetilde{N}+d+1}$ and $E_{\xi,\varepsilon} = E_{\xi,\varepsilon,\delta_\varepsilon}$. We obtain that $ E_{\xi,\varepsilon} \in C^{N}(\overline{X_n})$ with $P(D) E_{\xi,\varepsilon} = \delta_\xi$ in $X_k$. Furthermore, there are $C_1,C_2 >0$ such that for all $\varepsilon \in (0,\varepsilon_0)$
	$$
	\|E_{\xi,\varepsilon,\delta}\|_{n,N} \leq C_1 \varepsilon \qquad \mbox{and} \qquad \|E_{\xi,\varepsilon,\delta}\|_{k,K}^* \leq \frac{C_2}{\varepsilon^s}
	$$
	with $s = \widetilde{s}(N+\widetilde{N}+d+1)$. This completes the proof. \hfill \qed
	
	Theorem \ref{thm: main2}(c)  can be shown in the same way as  Theorem \ref{thm: main2}(b) (see in particular the values of $s$ in terms of $\widetilde{s}$ in the above proof). We leave the details to the reader.

	\section{The condition $(\Omega)$ for $\mathscr{E}_P(X)$ if $X$ is convex}\label{sect:convex}
	In this final section we use Theorem \ref{thm: main2}(a) to prove that $P(D): \mathscr{E}(X) \to  \mathscr{E}(X)$ is surjective and $\mathscr{E}_P(X)$ satisfies $(\Omega)$ for any non-zero differential operator $P(D)$ and any open convex set $X \subseteq \R^d$. To this end, we show that \eqref{eq: Omega} holds  for any exhaustion by relatively compact open convex subsets $(X_n)_{n\in\N}$ of $X$. The latter is a consequence of the following lemma.
	
	\begin{lemma}\label{FS}
		Let $P\in\C[\xi_1,\ldots,\xi_d] \backslash \{0\}$. Let $K\subseteq \R^d$ be  compact  and convex,  and let $\xi \in \R^d$  be such that $\xi \notin K$.  For all $\varepsilon \in (0,1)$ there exists $E_{\xi,\varepsilon} \in  \mathscr{D}'(\R^d)$ with $P(D)E_{\xi,\varepsilon} = \delta_\xi$ in $\R^d$  such that 
		$$
		\| E_{\xi,\varepsilon} \|^*_{K,d+1} \leq \varepsilon 
		$$
		and for every  $L \subset \R^d$ compact and convex there are $s,C > 0$ such that 
		$$
		\qquad \| E_{\xi,\varepsilon} \|^*_{L,d+1} \leq \frac{C}{\varepsilon^s}. 
		$$
	\end{lemma}
	The rest of this section is devoted to the proof of  Lemma \ref{FS}, which is based on a construction of fundamental solutions due to H\"ormander \cite[proof of Theorem 7.3.10]{HoermanderPDO1}. We need some preparation.  For $Q \in \C[\xi_1,\ldots,\xi_d]$ we define
	$$
	\widetilde{Q}(\zeta) = \left(\sum_{\alpha \in \N^d} |Q^{(\alpha)}(\zeta)|^2\right)^{1/2}, \qquad \zeta \in \C^d.
	$$
	Let $m \in \N$.  We denote by $\operatorname{Pol}^\circ(m)$ the finite-dimensional vector space of non-zero polynomials in $d$ variables of degree at most $m$ with the origin removed. By \cite[Lemma 7.3.11 and Lemma 7.3.12]{HoermanderPDO1} there exists a non-negative $\Phi \in C^\infty(\operatorname{Pol}^\circ(m) \times \C^d)$ such that
	\begin{itemize}
		\item[$(i)$] For all  $Q \in \operatorname{Pol}^\circ(m)$ it holds that $\Phi(Q, \zeta) = 0$ if $|\zeta| > 1$ and
		$$
		\int_{\C^d}\Phi(Q, \zeta) d\zeta = 1.
		$$
		\item[$(ii)$] For all entire functions $F$ on $\C^d$ and $Q \in \operatorname{Pol}^\circ(m)$ it holds that
		$$
		\int_{\C^d} F(\zeta) \Phi(Q, \zeta)d\zeta = F(0).
		$$
		\item[$(iii)$] There is $A >0$ such that for all $Q \in \operatorname{Pol}^\circ(m)$ and $\zeta \in \C^d$ with $\Phi(Q,\zeta) \neq 0$ it holds that
		$$
		\widetilde{Q}(0) \leq A|Q(\zeta)|.
		$$
	\end{itemize}
	
	Let $K \subseteq \R^d$ be compact and convex. As customary, we define the supporting function of $K$ as
	$$
	H_K(\eta) = \sup_{x \in K} \eta \cdot x, \qquad \eta \in \R^d.
	$$
	Note that $H_K$ is subadditive and positive homogeneous of degree $1$. Furthermore, it holds that \cite[Theorem 4.3.2]{HoermanderPDO1}
	\begin{equation}
		\label{support-convex}
		K = \{ x \in \R^d \, | \, \eta\cdot x \leq H_K(\eta), \, \forall \eta \in \R^d \}.
	\end{equation} 
	
	We define the Fourier transform of $\varphi \in \mathscr{D}(\R^d)$ as
	$$
	\widehat{\varphi}(\zeta) = \int_{\R^d} \varphi(x) e^{- i \zeta \cdot x} dx, \qquad \zeta \in \C^d.
	$$ 
	Then, $\widehat{\varphi}$ is an entire function on $\C^d$.  For all $N \in \N$ there is $C >0$ such that for all $\varphi \in \mathscr{D}(\R^d)$
	\begin{equation}
		\label{PWS}
		|\widehat{\varphi}(\zeta)| \leq C\|\varphi\|_{L^1,N} \frac{e^{H_{\operatorname{ch} \supp \varphi}( \operatorname{Im \zeta})}}{(2+|\zeta|)^{N}}, \qquad \zeta \in \C^d,
	\end{equation}
	where $\operatorname{ch} \supp \varphi$ denotes the convex hull of $\supp \varphi$. We are ready to show Lemma \ref{FS}.
	\begin{proof}
		We may assume without loss of generality that $\xi = 0$. Since $0 \notin K$, \eqref{support-convex} implies that there is $\eta \in \R^d$ such that  $H_K(-\eta) < 0$. For $t > 0$ and $\sigma \in \R^d$ we define $P_{t,\sigma} = P(\sigma  + it\eta + \, \cdot \,) \in \operatorname{Pol}^\circ(m)$. Note that there is $c >0$ such for all  $t > 0$ and $\sigma \in \R^d$
		\begin{equation}
			\label{lowerb}
			\widetilde{P_{t,\sigma}}(0) = \widetilde{P}(\sigma  + it\eta) \geq c.
		\end{equation}
		Let $\Phi$ be as above. We define $F_t \in \mathscr{D}'(\R^d)$ via (cf.\ \cite[proof of Theorem 7.3.10]{HoermanderPDO1})
		$$
		\langle  F_t, \varphi \rangle = \frac{1}{(2\pi)^d} \int_{\R^d} \int_{\C^d} \frac{\widehat{\varphi}(-\sigma - it\eta - \zeta)}{P(\sigma + it\eta + \zeta)} \Phi(P_{t,\sigma}, \zeta) d\zeta d\sigma, \qquad \varphi \in \mathscr{D}(\R^d).
		$$
		Let $L$ be an arbitrary compact convex subset of $\R^d$. By properties $(i)$ and $(iii)$ of $\Phi$, \eqref{PWS} (with $N = d+1$) and \eqref{lowerb} we have that for all $\varphi \in \mathscr{D}_L$
		\begin{eqnarray}
			\label{ess-bound}
			\nonumber |\langle  F_t, \varphi \rangle | &\leq&  \frac{1}{(2\pi)^d} \int_{\R^d} \int_{|\zeta| \leq 1} \frac{|\widehat{\varphi}(-\sigma - it\eta - \zeta)|}{|P_{t,\sigma}(\zeta)|} \Phi(P_{t,\sigma}, \zeta) d\zeta d\sigma \\ \nonumber
			&\leq&  \frac{AC\|\varphi\|_{L^1,d+1}}{(2\pi)^d} \int_{\R^d} \int_{|\zeta| \leq 1} \frac{e^{H_L(-t\eta - \operatorname{Im} \zeta)}}{(2+|\sigma + it\eta + \zeta|)^{d+1}\widetilde{P_{t,\sigma}}(0) }\Phi(P_{t,\sigma}, \zeta) d\zeta d\sigma \\ \nonumber
			&\leq&  \frac{AC\|\varphi\|_{L^1,d+1}}{(2\pi)^dc} \int_{\R^d} \int_{|\zeta| \leq 1} \frac{e^{tH_L(-\eta)}  e^{H_L(- \operatorname{Im} \zeta)}}{(1+|\sigma|)^{d+1}}\Phi(P_{t,\sigma}, \zeta) d\zeta d\sigma \\ 
			&\leq&  C'_L\|\varphi\|_{L^\infty,d+1} e^{tH_L(-\eta)},		    
		\end{eqnarray}
		where
		$$
		C'_L = \frac{AC|L|}{(2\pi)^dc} \max_{|\zeta| \leq 1 }e^{H_L(- \operatorname{Im} \zeta)} \int_{\R^d} \frac{1}{(1+|\sigma|)^{d+1}} d\sigma. 
		$$
		In particular, $F_t$ is a well-defined distribution. Property $(ii)$ of $\Phi$ and Cauchy's integral formula yield that for all $\varphi \in \mathscr{D}(\R^d)$
		\begin{eqnarray*}
			\langle P(D) F_t, \varphi \rangle  &=& \langle  F_t, P(-D)\varphi \rangle  \\	    
			&=& \frac{1}{(2\pi)^d} \int_{\R^d} \int_{\C^d} \widehat{\varphi}(-\sigma - it\eta - \zeta)\Phi(P_{t,\sigma}, \zeta) d\zeta d\sigma \\
			&=& \frac{1}{(2\pi)^d} \int_{\R^d}  \widehat{\varphi}(-\sigma - it\eta) d\sigma \\
			&=& \frac{1}{(2\pi)^d} \int_{\R^d}  \widehat{\varphi}(\sigma) d\sigma = \varphi(0)
		\end{eqnarray*}
		and thus $P(D) F_t = \delta$. For $\varepsilon \in (0,1)$ we set $t_\varepsilon = \log \varepsilon/ H_K(-\eta) > 0$ and $E_{0,\varepsilon} = F_{t_\varepsilon}$. Then, $P(D) E_{0,\varepsilon} = \delta$. By \eqref{ess-bound}, we obtain that for all $\varepsilon \in (0,1)$
		$$
		\| E_{0,\varepsilon} \|^*_{K,d+1} \leq  C'_K \varepsilon
		$$
		and, for any $L \subseteq \R^d$ compact and convex,
		$$
		\| E_{0,\varepsilon} \|^*_{L,d+1} \leq  \frac{C'_L}{\varepsilon^s},
		$$
		with $s = |H_L(-\eta)/ H_K(-\eta)|$.  This gives the desired result.
	\end{proof}


\begin{thebibliography}{999}
		\setlength{\itemsep}{0pt}
		
		\bibitem{B-D2006}J.~Bonet, P.~Doma\'nski, \emph{Parameter dependence of solutions of differential equations on spaces of distributions and the splitting of short exact sequences}, J. Funct. Anal. \textbf{230} (2006), 329--381.
		
		\bibitem{B-D2008}J.~Bonet, P.~Doma\'nski, \emph{The splitting of exact sequences of PLS-spaces and smooth dependence of solutions of linear partial differential equations}, Adv. Math. \textbf{217} (2008), 561--585. 
		
		\bibitem{DeKa22} A.~Debrouwere, T.~Kalmes, \emph{Linear topological invariants for kernels of convolution and differential operators}, J.\ Funct.\ Anal.\ \textbf{284} (2023), Paper No.\ 109886.
		
		\bibitem{DeKa22-2} A.~Debrouwere, T.~Kalmes, \emph{Quantitative Runge type approximation theorems for zero solutions of certain partial differential operators}, accepted for publication in Israel J.\ Math.
		
		\bibitem{Domanski} P.~Doma\'nski, \emph{Real analytic parameter dependence of solutions of differential equations}, Rev. Mat. Iberoam. \textbf{26} (2010), 175--238.
		
		\bibitem{Domanski2011} P.~Doma\'nski, \emph{Real analytic parameter dependence of solutions of differential equations over Roumieu classes},  Funct. et Approx. Comment. Math.  \textbf{26} (2011), 79--109.
		
		\bibitem{EnPe21}
		A.~Enciso, D.~Peralta-Salas, \emph{Approximation {T}heorems for the {S}chr\"{o}dinger {E}quation and {Q}uantum {V}ortex {R}econnection.} Comm.\ Math.\ Phys.\ \textbf{387} (2021), 1111--1149.
		
		\bibitem{FW04b} L.~Frerick, J.~Wengenroth, \emph{Partial differential operators modulo smooth functions}, Bull. Soc. Roy. Sci. Li\`ege \textbf{73} (2004), 119--127.
		
		\bibitem{FrKa}
		L.~Frerick, T.~Kalmes, \emph{Some results on surjectivity of augmented semi-elliptic differential operators}, Math. Ann.\ \textbf{347} (2010), 81--94.
		
		\bibitem{HoermanderPDO1} L.~H{\"o}rmander, \emph{The Analysis of Linear Partial Differential Operators, {I}}, Springer-Verlag, Berlin, 2003.
		
		\bibitem{HoermanderPDO2}
		L.~H{\"o}rmander, \emph{The Analysis of Linear Partial Differential Operators, {II}},  Springer-Verlag, Berlin, 2005.
		
		\bibitem{Kalmes12-2}
		T.~Kalmes, \emph{The augmented operator of a surjective partial differential operator with constant coefficients need not be surjective}, Bull.\ Lond.\ Math.\ Soc.\ \textbf{44} (2012), 610--614.
		
		\bibitem{Langenbruch1996} M.~Langenbruch, \emph{Surjective partial differential operators on spaces of ultradifferentiable functions of Roumieu type}, Results in Math. \textbf{29} (1996), 254--275.
		
		\bibitem{Langenbruch1998} M.~Langenbruch, \emph{Surjectivity of partial differential operators on Gevrey classes
			and extension of regularity}, Math. Nachr. \textbf{196} (1998), 103-140.
		
		\bibitem{Langenbruch2004} M.~Langenbruch, \emph{Characterization of surjective partial differential operators on spaces of real
			analytic functions}, Studia Math. \textbf{162} (2004), 53--96.
		
		\bibitem{MTV90} R.~Meise, B.~A.~Taylor, D.~Vogt, \emph{Characterization of the linear partial differential operators with constant coefficients that admit a continuous linear right
			inverse}, Ann. Inst. Fourier (Grenoble) \textbf{40} (1990), 619--655.
		
		\bibitem{MTV96} R.~Meise, B.~A.~Taylor, D.~Vogt, \emph{Continuous linear right inverses for partial differential operators on non-quasianalytic classes and on ultradistributions}, Math. Nachr. \textbf{196} (1998), 213-242.
		\bibitem{M-V} R.~Meise, D.~Vogt, \emph{Introduction to Functional Analysis}, Clarendon Press, Oxford, 1997.
		
		\bibitem{Petzsche1980}
		H.J.~Petzsche, \emph{Some results of Mittag-Leffler-type for vector valued functions and spaces of class $A$}, in: K.D.~Bierstedt, B.~Fuchssteiner (Eds.), \emph{Functional Analysis: Surveys and Recent Results}, North-Holland, Amsterdam, 1980, pp.\ 183--204.
		
		\bibitem{RuSa19}
		A.~R\"{u}land, M.~ Salo, \emph{Quantitative Runge approximation and inverse problems}, Int.\ Math.\ Res.\ Not.\ IMRN \textbf{20} (2019), 6216--6234.
		
		\bibitem{RuSa20}
		A.~R\"{u}land, M.~ Salo, \emph{The fractional Calder\'on problem: Low regularity and stability}, Nonlinear Anal. \textbf{93} (2020), 111529.
		
		\bibitem{Schwartz} L.~Schwartz, \emph{Th\'eorie des Distributions}, Hermann, Paris, 1966.
		
		\bibitem{Vogt1983-2} D.~Vogt, \emph{On the solvability of $P(D) f = g$ for vector valued functions}, RIMS Kokyoroku \textbf{508} (1983), 168--181.
		
		\bibitem{Vogt1987} D.~Vogt, \emph{On the functors $\operatorname{Ext}^{1}(E, F)$ for Fr\'{e}chet spaces},  Studia Math. \textbf{85} (1987), 163--197.
		
		\bibitem{Vogt2006} D.~Vogt, \emph{Invariants and spaces of zero solutions of linear partial differential operators}, Arch. Math. \textbf{87} (2006), 163--171.
		
		\bibitem{Wengenroth03} J.~Wengenroth, \emph{Derived Functors in Functional Analysis}, Springer-Verlag, Berlin, 2003.
		
		\bibitem{Wengenroth11} J.~Wengenroth, \emph{Surjectivity of partial differential operators with good fundamental solutions}, J. Math. Anal. Appl. \textbf{379} (2011), 719--723.
	\end{thebibliography}
\end{document}